\documentclass[11pt,letterpaper]{article}
\makeatletter
\usepackage[dvips]{graphicx}
\usepackage{latexsym,amssymb,amsmath,amsfonts,amsthm}%amssymb,vatola, vatola,amsthm,amssymb

\newtheorem{thm}{\textbf Theorem}[section]
\newtheorem{lem}{\textbf Lemma}[section]
\newtheorem{rem}{\textbf Remark}[section]

\newtheorem{prop}{\textbf Proposition}[section]
\newtheorem{defin}{\textbf Definition}[section]

\newcommand{\md}{\mbox{d}}
\newcommand{\be}{\begin{eqnarray}}
\newcommand{\ee}{\end{eqnarray}}
\newcommand{\mr}{\mathbb{R}}
\newcommand{\mc}{\mathcal}
\newcommand{\mb}{\mathbb}
\newcommand{\mx}{\mbox}
\newcommand{\ep}{\varepsilon}
\newcommand{\tr}{\triangle}
\newcommand{\pt}{\partial}
\newcommand{\bes}{\begin{eqnarray*}}
\newcommand{\ees}{\end{eqnarray*}}
\newcommand{\Om}{\Omega}
\newcommand{\om}{\omega}

\newcommand{\al}{\alpha}
\newcommand{\la}{\lambda}
\newcommand{\La}{\Lambda}
\newcommand{\De}{\Delta}
\newcommand{\ga}{\gamma}
\newcommand{\de}{\delta}

\begin{document}
\begin{titlepage}
\title{\bf Global well-posedness of incompressible flow in porous media with critical diffusion in Besov spaces}
\author{Baoquan Yuan$^{1,2}$ and Jia Yuan$^3$
\\ $^1$Institute of Applied Physics and Computational Mathematics,\\ Beijing 100088,
        China\\ $^2$School of Mathematics and Informatics,
     \\  Henan Polytechnic University, Henan, 454000, China.\\
$^3$Department of  Mathematics, Beihang University,
\\ LMIB of the Ministry of Education; Beijing 100083, China.\\
         (bqyuan@hpu.edu.cn)}

\date{}
\end{titlepage}
\maketitle
\begin{abstract}
In this paper we study the model of heat transfer in a porous medium
with a critical diffusion. We obtain global existence and uniqueness
of solutions to the equations of heat transfer of incompressible
fluid in Besov spaces $\dot B^{3/p}_{p,1}(\mr^3)$ with $1\le
p\le\infty$ by the method of modulus of continuity and Fourier
localization technique.
 \vskip0.1in

\noindent{\bf AMS Subject Classification 2000:}\quad 76S05, 76D03
\end{abstract}

\vspace{.2in} {\bf Key words:}\quad Flows in porous media, global
well-posedness, Besov spaces.

%blow-up criterion, regularity of weak solutions.

%%\newpage

%%%%%%%%%%%%%%%%%%%%%%%%%%%%%%%%%%%%%%%%%%%%%%%%%%%%%%%%%%%%%%%%%%%%%%%%%%%%%%%%%%%%%%%%%%%%%%%%%%%%%%%%%%%%%%%%%%%%%%%%%%%%

\section{Introduction}
\setcounter{equation}{0}

In this paper we consider the transfer of the heat with a general
diffusion term in an incompressible flow in the porous medium. The
equations is the following:
 \be\label{PE}
\begin{cases}
\frac{\pt\theta}{\pt t}+u\cdot\nabla\theta+\nu\La^\al\theta=0,\
x\in
\mr^3,\ t>0,\\
u=-k(\nabla p+g\ga\theta),x\in
\mr^3,\ t>0,\\
\mx{div}u=0,\\
\theta(0,x)=\theta_0.
\end{cases}
 \ee
Here $\nu>0$ is the dissipative coefficient and $\ga$ is the
matrix medium permeability in the different directions,
respectively, divided by the viscosity, $g$ is the acceleration
due to gravity and the vector $\ga\in \mr^3$ is the last canonical
vector $e_3$, $\theta$ is the liquid temperature and $u=-k(\nabla
p+g\ga\theta)$ the liquid discharge by the Darcy's law, $p$ is the
pressure of the liquid. For more details see \cite{N-B}. To
simplify the notation, we set $k=g=1$.

  The operator $\La^\al=(-\De)^{\al/2}$ is defined by the Fourier transform
  \be
\mc{F}(\La^\al \theta)(\xi)=|\xi|^\al\mc{F}\theta(\xi),
  \ee
for $0\le \al\le 2$.

The case $\al=1$ is called the critical case, $1<\al\le 2$ is the
sub-critical case and $0\le \al<1$ is the super-critical case.

According to the Darcy's law and the incompressibility condition,
one has
 \bes
\De u=-\mx{curl}(\mx{curl} u)=\bigg(\frac{\pt^2\theta}{\pt x_1\pt
x_3},\frac{\pt^2\theta}{\pt x_2\pt x_3},-\frac{\pt^2\theta}{\pt
x_1^2}-\frac{\pt^2\theta}{\pt x^2_2}\bigg).
 \ees

By Newton potential formula and integrating by parts one has
 \be\nonumber
u&=&-\frac23(0,0,\theta)+\frac1{4\pi}\mx{P.V.}\int_{\mr^3}K(x-y)\theta(t,y)\md
y,\ x\in\mr^3,\\ \label{1.3}
&=&c(\theta)+\mc{P}(\theta),
 \ee
where
 \bes
K(x)=\bigg(\frac{3x_1x_3}{|x|^5},\frac{3x_2x_3}{|x|^5},\frac{2x_3^2-x_1^2-x_2^2}{|x|^5}\bigg)
 \ees
is the kernel function of a singular integral $\mc{P}(\theta)$,
for detail see \cite{C-C-G-O}.

When $\nu=0$ and space dimension $n=2$, D. C\'{o}rdoba and F.
Gancedo \cite{C-G-O} obtained the local existence and uniqueness
by the particle trajectory method in H\"{o}lder space$C^s$ for
$0<s<1$ and gave some blow-up criteria of smooth solutions, for
example, the blow-up criterion in BMO space similar to the Euler
equations and the geometric constraint conditions under which no
singularities are possible, for details see \cite{C-G-O} and
reference therein.

When $\nu>0$, A. Castro, D. C\'{o}rdoba, F. Gancedo and R. Orive
\cite{C-C-G-O} constructed the global solutions to (\ref{PE}) in
the Sobolev space $H^s$ with $s>0$ for the sub-critical diffusion
case. In the super-critical diffusion case the global
well-posedness for small initial data in $H^s$ with $s>n/2+1$ and
the local well-posedness in the space $H^s$ with $s>(n-\al)/2+1$
were obtained in \cite{C-C-G-O}. In the critical diffusion case
the global well-posedness of smooth solutions can also be obtained
for smooth initial data by the method as in \cite{K-N-V,C-V} for
the critical dissipative quasi-geostrophic equations.

Before presenting our method and results let us first clarify the
notion of critical space (super-critical and sub-critical spaces,
respectively). If $\theta(t,x)$ is a solution to the equations
(\ref{PE}), then $\theta_\la=\la^{\al-1}\theta(\la^\al t,\la x)$ is
also a solution for $\la>0$. A translation invariant homogeneous
Banach space of distribution $X$ is called a critical space, if its
norm is invariant under the scaling transform
$f_\la=\la^{\al-1}f(\la x )$, i.e. $\|f_\la\|_X=\|f\|_X$ for any
$\la>0$. Similarly, it is called a super-critical space
(sub-critical space), if $\log_{\la}\frac{\|f_\la\|_X}{\|f\|_X}<0\
(>0)$ for $\la>0$. Noting that the space $H^s$ for $s>(n-\al)/2+1$
is a sub-critical space for the super-critical or critical diffusion
cases of the incompressible flow equations in the porous medium, so
the energy method and Sobolev estimates are available. But for the
critical space $\dot B^{3/p}_{p,1}(\mr^3)$ it needs a different
method. Indeed, the energy methods are not applicable. One needs
first establish the local existence, uniqueness and higher
regularity based on a priori estimate of the following
transport-diffusion equation
 \bes
\begin{cases}
\pt_tu+v\cdot\nabla u+\nu\La^\al u=f,\\
u(0,x)=u_0.
\end{cases}
 \ees
That is
 \bes
\nu^{1/r}\|u\|_{\tilde L^r_T\dot B^{s+\al/r}_{p,q}}\le
C\mx{e}^{CZ(T)}\Big(\|u_0\|_{\dot
B^s_{p,q}}+\nu^{1/r_1-1}\|f\|_{\tilde L^{r_1}_T\dot
B^{s-\al+\al/r_1}_{p,q}}\Big),
 \ees
where $Z(T)=\int^T_0\|\nabla v(t)\|_{\dot
B^{n/p_1}_{p_1,\infty}\cap L^\infty}\md t$. For details see
Proposition \ref{prop2}.

By virtue of the method of modulus of continuity \cite{K-N-V}, we
prove the global existence and uniqueness of solutions to the
equations (\ref{PE}) with $\al=1$ in critical Besov space $\dot
B^{3/p}_{p,1}(\mr^3)$ with $1\le p\le\infty$. The key point is to
construct a new modulus of continuity, which control the blow-up
of the smooth local solutions to the equations (\ref{PE}). Assume
that $\theta$ has a modulus $\om$, Kiselev, Nazarov and Volberg in
\cite{K-N-V} proved that the Riesz transform $R_j(\theta)$ had a
modulus of continuity
 \be\label{RModulus}
\Om_1(x)=A\bigg(\int^x_0\frac{\om(s)}s\md
s+x\int^\infty_x\frac{\om(s)}{s^2}\md s\bigg),
 \ee
where $A$ is a constant.
 Noticing the
relation (\ref{1.3}) of $u$ and $\theta$ which is equivalent to
double Riesz transforms. We prove that the singular integral
operator $\mc{P}(\theta)$ in (\ref{1.3}) do not spoil modulus of
continuity of $\theta$ too much. In fact, it has a modulus of
continuity
 \be\label{Modulus}
\Om(x)=C\bigg(\int^x_0\frac{\om(s)}s\log\Big(\frac{\mx{e}x}
s\Big)\md
s+x\int^\infty_x\frac{\om(s)}{s^2}\log\Big(\frac{\mx{e}s}x\Big)\md
s\bigg),
 \ee
where $C>0$ is a constant. See Lemma \ref{lem3.1}.

Comparing with (\ref{RModulus}), the formula (\ref{Modulus}) of
modulus of continuity of a double Riesz transform has an additional
term $\log\frac{\mx{e}x}s$ or $\log\frac{\mx{e}s}x$, which requires
us to construct the modulus of continuity $\omega$ of temperature
$\theta$ to be slowly increased at infinite. In fact we construct
the continuous function $\om(x)$ as follows
 \be
\om(x)=\begin{cases} x-x^{3/2},\ \ &\mx{ if }0\le x\le \de,\\
\de-\de^{3/2}+\frac\ga 3\arctan\frac{1+\log \frac x\de}3-\frac\ga
3\arctan\frac13,\ \ &\mx{ if } \de\le x,
\end{cases}
 \ee
which is a increasing bounded concave function when
$x\rightarrow\infty$. While the modulus of continuity of $\theta$ in
\cite{K-N-V} has a double logarithm-type increase at infinite
 \bes
\om_1(x)=\begin{cases} x-x^{3/2},\ \ &\mx{ if }0\le x\le \de,\\
\de-\de^{3/2}+\ga\log\bigg(1+\frac14\log\frac x\de\bigg),\ \ &\mx{
if } \de\le x.
\end{cases}
 \ees

%\begin{figure}[htpb]
%\mbox{}\hskip12mm\epsfig{figure=photo.eps} \caption{Illustration
%for the definitions}
%\end{figure}

%\includegraphics[scale=1]{photo.eps}
%Pictures of functions $\om_1(x)$ and $\om(x)$, $1\le x\le 30 $
%
To this end we present our main result in the following Theorem
\ref{thm1.1}.

\begin{thm}\label{thm1.1}
Let $\theta_0\in \dot B^{3/p}_{p,1}(\mr^3)$ with $1\le p\le
\infty$, then the critical diffusion equations (\ref{PE}) of heat
transfer of incompressible fluid possesses a unique global
solution $\theta\in C(\mr^+;\dot B^{3/p}_{p,1}(\mr^3))\cap
L^1_{loc}(\mr^+;\dot B^{3/p+1}_{p,1}(\mr^3))$.
\end{thm}

To prove our main theorem, we need a local existence theorem as
follows.

\begin{thm}\label{thm1.2}
Let $\theta_0\in \dot B^{3/p}_{p,1}(\mr^3)$ with $1\le p\le \infty$,
then there exists a time $T>0$ such that the equations (\ref{PE})
possesses a unique local solution $\theta\in C([0,T);\dot
B^{3/p}_{p,1})$ satisfying
 \be
\theta\in \tilde L^\infty([0,T);\dot B^{3/p}_{p,1}(\mr^3))\cap
L^1((0,T);\dot B^{3/p+1}_{p,1}(\mr^3)).
 \ee
Furthermore, we also have $t^\beta\theta\in \tilde
L^\infty((0,T);\dot B^{3/p+\beta}_{p,1}(\mr^3))$ for $\beta> 0$.
\end{thm}

For the initial data $\theta_0$ which satisfies the condition
$\|\theta_0\|_{\dot B^{3/p}_{p,1}}<\infty$, global well-posedness
can not be obtained. So we give a blow-up criterion of smooth
solutions in the following Theorem \ref{thm1.3}.

\begin{thm}\label{thm1.3}
Let $T>0$ and $\theta_0\in \dot B^{3/p}_{p,1}(\mr^3)$ with $1\le
p\le \infty$. Assume that $\theta \in \tilde L^\infty([0,T);\dot
B^{3/p}_{p,1}(\mr^3))\cap L^1((0,T);\dot B^{3/p+1}_{p,1}(\mr^3))$ is
a smooth solution to the equations (\ref{PE}), if $\theta$ satisfies
 \be\label{1.5}
\int^T_0\|\nabla\theta(t)\|_{L^\infty}\md t<\infty,
 \ee
then $\theta(t,x)$ can be continually extended to the interval
$(0,T')$ for some $T'>T$.
\end{thm}

\begin{rem}\label{rem1.1}
Actually, we have more general blow-up criterion. But in our case,
the $L^\infty$ norm is enough. Indeed, assume that $\theta\in
C([0,T];H^s(\mr^3))$ with $s>\frac n2+1$, if $\theta$ satisfies
 \be
\int^T_0\|\nabla\theta(t)\|_{\dot B^0_{\infty,\infty}}\md t<\infty,
 \ee
then $\theta(t,x)$ can be continually extended to the interval
$(0,T')$ for some $T'>T$.

The proof is standard. Energy method and the following logarithmic
Sobolev inequality \cite{Yuan, K-O-T}
 \bes
\|f\|_{L^\infty}\le C(1+\|f\|_{\dot
B^0_{\infty,\infty}}\log(\mx{e}+\|f\|_{H^s}))
 \ees
 for $s>\frac n2$, immediately yield the result.
\end{rem}

%%%%%%%%%%%%%%%%%%%%%%%%%%%%%%%%%%%%%%%%%%%%%%%%%%%%%%%%%%%%%% Section 2 %%%%%%%%%%%%%%%%%%%%%%%%%%%%%%%%%%%%%%%%%%%%%%%%

\section{Preliminaries}
\setcounter{equation}{0}

We first introduce the Littlewood-Paley decomposition and
definition of Besov spaces. Given $f\in \mc{S}(\mr^n)$ the
Schwartz class of rapidly decreasing function, define the Fourier
transform as
 \bes
\hat{f}(\xi)=\mc{F}f(\xi)=(2\pi)^{-n/2}\int_{\mr^n}\mx{e}^{-ix\cdot\xi}f(x)\md
x,
 \ees
and its inverse Fourier transform:
 \bes
\check{f}(x)=\mc{F}^{-1}f(x)=(2\pi)^{-n/2}\int_{\mr^n}\mx{e}^{ix\cdot\xi}f(\xi)\md\xi.
 \ees
Choose a nonnegative radial function $\chi\in C^\infty_0(\mr^n)$
such that $0\le \chi(\xi) \le 1$ and
 \bes
\chi(\xi)=
\begin{cases}
1, &\mx{ for } |\xi|\le \frac34,\\
0, &\mx{ for } |\xi|>\frac 43,
\end{cases}
 \ees
and let $\hat{\varphi}(\xi)=\chi(\xi/2)-\chi(\xi)$,
$\chi_j(\xi)=\chi(\frac{\xi}{2^j})$ and
$\hat{\varphi}_j(\xi)=\hat{\varphi}(\frac \xi{2^j})$ for $j\in
\mb{Z}$. Write
 \bes
h(x)&=&\mc{F}^{-1}\chi(x),\ h_j(x)=2^{nj}h(2^jx);\\
\varphi_j(x)&=&2^{nj}\varphi(2^jx).
 \ees
Define the Littlewood-Paley operators $S_j$ and $\tr_j$,
respectively, as
 \bes
 \tr_{-1}u(x)&=&h*u(x),\\
\tr_ju(x)&=&\varphi_j*u(x)=S_{j+1}u(x)-S_ju(x), \ \mx{ for } j\ge
0,\\
\tr_ju(x)&=&0,\ \mx{ for } j\le -2,\\
 S_ju(x)&=&(1-\sum_{k\ge j}\tr_k)u(x),\ \mx{ for } j\in \mb{Z}.
 \ees
Formally $\tr_j$ is a frequency projection to the annulus
$|\xi|\approx 2^j$, while $S_j$ is a frequency projection to the
ball $|\xi|\lesssim 2^j$ for $j\in \mb{Z}$. For any $u(x)\in
L^2(\mr^n)$ we have the Littlewood-Paley decomposition
 \be\nonumber
u(x)&=&h*u(x)+\sum_{j\ge 0}\varphi_j*u(x)\ (\mx{Inhomogeneous
decomposition}),\\
\label{homogeneous} u(x)&=&\sum^\infty_{j=-\infty}\varphi_j*u(x)\
(\mx{Homogeneous decomposition}).
 \ee
Here homogeneous decomposition (2.1) holds in the sense of modulus
of polynomial function. Clearly,
 $\mx{supp}\chi(\xi)\cap \mx{supp}\hat{\varphi}_j(\xi)=\emptyset$,
 for $j\ge 1$, $\mx{supp}\hat{\varphi}_j(\xi)\cap
\mx{supp}\hat{\varphi}_{j'}(\xi)=\emptyset$, for $|j-j'|\ge 2$,
and $\tr_j(S_{k-1}u\tr_ku)=0$ for $|j-k|\ge5$.

Next, we recall the definition of Besov spaces. Let $s\in \mr$ and
$1\le p,\ q\le +\infty$, the Besov space $B^s_{p,q}(\mr^n)$
abbreviated as $B^s_{p,q}$ is defined by
 \bes
B^s_{p,q}=\{f(x)\in \mc{S}'(\mr^n);\|f\|_{B^s_{p,q}}<+\infty \},
 \ees
where
 \bes
\|f\|_{B^s_{p,q}}=\begin{cases} \|h*f\|_p+\bigg(\sum_{j\ge
0}2^{jsq}\|\varphi_j*f\|^q_p\bigg)^{1/q}, \ \ &\mx{ for
}q<+\infty,\\\|h*f\|_p+ \sup_{j\ge 0}2^{js}\|\varphi_j*f\|_p, \ \
&\mx{ for }q=+\infty
\end{cases}
 \ees
is the Besov norm. The homogeneous Besov space $\dot{B}^s_{p,q}$
is defined by the dyadic decomposition as
 \bes
\dot{B}^s_{p,q}=\{f(x)\in \mathcal{Z}'(\mathbb{R}^n);\
\|f\|_{\dot{B}^s_{p,q}}<+\infty\},
 \ees
where
 \be\label{Besov}
\|f\|_{\dot{B}^s_{p,q}}=\begin{cases}
\bigg(\sum^\infty_{j=-\infty}2^{jsq}\|\varphi_j*f\|^q_p\bigg)^{1/q},\
\ &\mx{
for }q<+\infty,\\
\sup_{j\in\mathbb{Z}}2^{js}\|\varphi_j*f\|_p,\ \ &\mx{ for
}q=+\infty
\end{cases}
 \ee
is the homogeneous Besov norm, and $\mathcal{Z}'(\mathbb{R}^n)$
denotes the dual space of $\mathcal{Z}(\mathbb{R}^n)=\{f(x)\in
\mathcal{S}(\mathbb{R}^n);\ D^\alpha \hat{f}(0)=0,\ \mbox{for any}
\alpha \in \mathbb{N}^n\ \mbox{multi-index}\}$ and can be identified
by the quotient space $\mathcal{S}'/\mathcal{P}$ with the polynomial
functional set $\mathcal{P}$. For details see \cite{Miao} and
\cite{Tr}.

\begin{rem}
The above definition does not depend on the choice of the radial
function $\chi$, and $\dot B^s_{p,q}(\mr^n)$ is a Banach space if
$s<\frac np$ or $s=\frac np$ and $q=1$.
\end{rem}

For the convenience, we recall the definition of Bony's para-product
formula which gives the decomposition of the product $f g$ of two
distributions $f$ and $g$.

%%%%%%%%%%%%%%%%%%%%%%%%Definition2.1

\begin{defin}
The para-product of two distributions $f$ and $g$ is defined by
 \bes
T_gf=\sum_{i\le j-2}\triangle_ig\triangle_jf=\sum_{j\in
\mathbb{Z}}S_{j-1}g\triangle_jf.
 \ees
The remainder of the para-product is defined by
 \bes
R(f,g)=\sum_{|i-j|\le1}\triangle_ig\triangle_jf.
 \ees
Then Bony's para-product formula reads
\begin{eqnarray}\label{2.4}
f g=T_gf+T_fg+R(f,g).
\end{eqnarray}
\end{defin}

Next we define two kinds of space-time Besov spaces that will be
used in our studies.

\begin{defin}
(1) Let $T>0,\ s\in\mr$ and $1\le p,\ q,\ r\le\infty$, $u(t,x)\in
\mc{S}'(\mr^4)$. We call $u(t,x)\in L^r(0,T;\dot
B^s_{p,q}(\mr^3))$ if and only if
 \be
\|u\|_{L^r\dot B^s_{p,q}}\triangleq \bigg\|\Big(\sum _{j\in
\mb{Z}}2^{jsq}\|\tr_ju\|_{L^p}^q\Big)^{1/q}\bigg\|_{L^r_T}<\infty.
 \ee
(2) Let $T,\ s,\ p,\ q,\ r$ and $u(t,x)$ be as in (1), we call
$u(t,x)\in \tilde L^r(0,T;\dot B^s_{p,q}(\mr^3))$ if and only if
 \be
 \|u\|_{\tilde L^r\dot B^s_{p,q}}\triangleq
 \Big(\sum_{j\in
 \mb{Z}}2^{jsq}\|\tr_ju\|^q_{L^r_TL^p}\Big)^{1/q}<\infty.
 \ee
\end{defin}

Obviously, by the Minkowski inequality, we have the relations
between the above two kinds of mixed space-time Besov spaces:
 \be
\|u\|_{\tilde L^r\dot B^s_{p,q}}\le \|u\|_{L^r\dot B^s_{p,q}},\ \
if\ \ r\le q,
 \ee
and
 \be
\|u\|_{L^r\dot B^s_{p,q}}\le \|u\|_{\tilde L^r\dot B^s_{p,q}},\ \
if\ \ q\le r.
 \ee
 We now recall some properties of the Besov spaces, for
details see \cite{Miao} or \cite{Tr}.

\begin{prop}\label{prop1}
The following properties of the Besov spaces hold:

(1). Let $\al\in \mr$, then the operator $\La^\al$ is an
isomorphism from $\dot B^s_{p,q}$ to $\dot B^{s-\al}_{p,q}$.

(2). If $p_1\le p_2$ and $q_1\le q_2$, then $\dot
B^s_{p_1,q_1}\hookrightarrow \dot
B^{s-n(\frac1{p1}-\frac1{p_2})}_{p_2,q_2}$.

(3). If $1\le p,\ q\le \infty$, $s>0$, $\alpha>0$ and $\beta>0$, and
$1\le p_i,\ q_i\le \infty$ ($i=1,\ 2,\ 3,\ 4$) so that
 \be\nonumber
\frac1p=\frac1{p_1}+\frac1{p_2}=\frac1{p_3}+\frac1{p_4},\\
\frac1q=\frac1{q_1}+\frac1{q_2}=\frac1q_3+\frac1{q_4}.
 \ee
Then there exists a constant $C$ such that $f_1\cdot f_2\in
\dot{B}^s_{p,q}(\mr^n)$ and
 \be\label{product}
\|f_1 f_2\|_{\dot{B}^s_{p,q}}\le C
(\|f_1\|_{\dot{B}^{s+\alpha}_{p_1,q_1}}\|f_2\|_{\dot{B}^{-\alpha}_{p_2,q_2}}+
\|f_1\|_{\dot{B}^{-\beta}_{p_3,q_3}}\|f_2\|_{\dot{B}^{s+\beta}_{p_4,q_4}}),
 \ee
for any $f_1\in \dot{B}^{s+\alpha}_{p_1,q_1}\cap
\dot{B}^{-\beta}_{p_3,q_3}$, $f_2\in
\dot{B}^{s+\beta}_{p_4,q_4}\cap \dot{B}^{-\alpha}_{p_2,q_2}$. If
$\al=0$, $p_2=q_2=\infty$ and $\beta=0$, $p_3=q_3=\infty$, then we
also have
 \be\label{product1}
\|f_1 f_2\|_{\dot{B}^s_{p,q}}\le C
(\|f_1\|_{\dot{B}^s_{p,q}}\|f_2\|_{L^\infty}+
\|f_1\|_{L^\infty}\|f_2\|_{\dot{B}^s_{p,q}}).
 \ee
\end{prop}
Using the para-product decomposition (\ref{2.4}) one can easily
prove the equality (\ref{product})-(\ref{product1}), and for the
proof of equality (\ref{product}) see \cite{Y-Z1}

In the following Lemma \ref{lem2.1} we recall the Bernstein
inequality which will be frequently used.

\begin{lem}\label{lem2.1}
Let $f\in L^p(\mr^n)$ with $1\le p\le q\le +\infty$ and $0<r<R$.
Then there exists constants $C>0$ and $C_k>0$ such that for any
$k\in \mb{N}$ and $\la>0$, one has
 \be\label{2.2}
\sup_{|\beta|=k}\|\pt^\beta f\|_q\le C\la^{k+n(1/p-1/q)}\|f\|_p, \
\mx{ if } \mx{supp}\hat f\subseteq \{\xi:\ |\xi|\le \la r\};
 \ee
 \be\label{2.3}
C^{-1}_k\la^k\|f\|_p\le\sup_{|\beta|=k}\|\pt^\beta f\|_p\le
C_k\la^k\|f\|_p,\ \mx{ if }\mx{supp}\hat f\subseteq \{\xi:\ \la
r\le |\xi|\le \la R\}.
 \ee
\end{lem}

The following lemmas will be useful in our discussions.

\begin{lem}$^{\cite{H-K,W-Y}}$\ \label{lem2.2}
Let $\psi$ be a smooth function supported on the shell $\{x\in
\mr^3:\ R_1\le |x|\le R_2,\ 0<R_1<R_2\}$. Then there exist two
positive constants $\mu$ and $C$ depending only on $\psi$ so that
for all $1\le p\le\infty$, $\al>0$, $t>0$ and $\la>0$, one has
 \bes
\|\psi(\la^{-1}D)\mx{e}^{-t\La^\al}u\|_{L^p}\le C \mx{e}^{-\mu
t\la^\al}\|\psi(\la^{-1}D)u\|_{L^p}.
 \ees
\end{lem}

\begin{lem}$^{\cite{Danchin}}$\label{lem2.3}
Let $v$ be a smooth vector field, and $\phi$ a solution to the
ordinary differential equation
 \be\nonumber
\begin{cases}
\frac{\md \phi(t,x)}{\md t}=v(t,\phi(t,x)),\\
\phi(0,x)=x.
\end{cases}
 \ee
Then for all $0\le t<\infty$, the flow $\phi(t,x)$ is a
diffeomorphism over $\mr^3$ and the following estimates hold:
 \bes
\|\nabla\phi(t)^{\pm 1}\|_{L^\infty}\le \mx{e}^{V(t)},
 \ees
 \bes
\|\nabla\phi(t)^{\pm 1}-Id\|_{L^\infty}\le \mx{e}^{2V(t)}-1,
 \ees
 \bes
\|\nabla^2\phi(t)^{\pm
1}\|_{L^\infty}\le\mx{e}^{V(t)}\int^t_0\|\nabla^2v(s)\|_{L^\infty}\mx{e}^{V(s)}\md
s,
 \ees
 where $V(t)=\int^t_0\|\nabla v(s)\|_{L^\infty}\md s$.
\end{lem}

%%%%%%%%%%%%%%%%%%%%%%%%%%%%%%%Lemma

The next lemma shows a estimate of exchange between $\La^s$ and
the flow $\phi_j$.

\begin{lem}$^{\cite{C-M-W}}$\label{lem2.4}
Let $v$ be a given vector field belonging to $L^1_{loc}(\mr^+;Lip)$,
$u_j\triangleq\tr_ju$, and $\phi_j$ denote the flow of the
regularized vector field $S_jv$ for $j\in \mb{Z}$. Then, for $u\in
\dot B^s_{p,\infty}$ with $0\le s<2$ and $1\le p\le\infty$, it
holds:
 \bes
\|\La^s(u_j\circ \phi_j)-(\La^su_j)\circ\phi_j\|_{L^p}\le
C2^{js}\mx{e}^{cV(t)}V^{1-\frac s2}(t)\|u_j\|_{L^p},
 \ees
where $V(t)=\int^t_0\|\nabla v(\tau)\|_{L^\infty}\md \tau$ and the
constant $C=C(s,p)$ depends only on $s$ and $p$.
\end{lem}

For the proofs of Lemma \ref{lem2.2}, \ref{lem2.3} and
\ref{lem2.4}, see \cite{H-K,W-Y}, \cite{Danchin} and \cite{C-M-W},
respectively.

Before we present the local existence of solutions to the
equations (\ref{PE}), we recall an optimal a priori estimate for
the following transport-diffusion equations in $\mr^n$:
 \be \label{TD}
\begin{cases}
\frac{\pt u}{\pt t}+v\cdot\nabla u-\nu\La^\al u=f\\
u(0,x)=u_0(x),
\end{cases}
 \ee
where $v$ is a fixed vector field which does not need to be
divergence free, $u_0$ is the initial data, $f$ is a given
external force term, and $\nu>0$ is a dissipative coefficient for
$0\le\al\le 2$.

\begin{prop}\label{prop2}
Let $1\le r_1\le r\le\infty$, $1\le p\le p_1\le\infty$ and $1\le
q\le\infty$. Assume $s\in\mr$ satisfies the following conditions:
 \be\label{Condition}
\begin{cases}
s<1+\frac n{p_1},\ (\mx{or }s\le 1+\frac n{p_1}, \mx{ if }q=1),\\
s>-\mx{min}\Big(\frac n{p_1},\frac n{p'}\Big),\ (\mx{or
}s>-1-\mx{min}\Big(\frac n{p_1},\frac n{p'}\Big), \mx{ if }
\mx{div }v=0).
\end{cases}
 \ee
There exists a constant $C>0$ depending only on $n$, $\al$, $s$,
$p$, $p_1$ and $q$, such that for any smooth solution $u$ of
equation (\ref{TD}), the following a priori estimate holds:
 \be\label{TDE}
\nu^{1/r}\|u\|_{\tilde L^r_T\dot B^{s+\al/r}_{p,q}}\le
C\mx{e}^{CZ(T)}\Big(\|u_0\|_{\dot
B^s_{p,q}}+\nu^{1/r_1-1}\|f\|_{\tilde L^{r_1}_T\dot
B^{s-\al+\al/r_1}_{p,q}}\Big),
 \ee
where $Z(T)=\int^T_0\|\nabla v(t)\|_{\dot B^{n/p_1}_{p_1,\infty}\cap
L^\infty}\md t$.

Moreover, if $u=v$, then for all $s>0$ ($s>-1$ if div $v=0$), the
estimate (\ref{TDE}) holds with $Z(T)=\int^T_0\|\nabla
v(t)\|_{L^\infty}\md t$.
\end{prop}

\begin{rem}
Danchin proved Proposition \ref{prop2} in \cite{Danchin} for $\al=2$
in the inhomogeneous Besov spaces case. Miao and Wu proved
Proposition \ref{prop2} in \cite{M-W} for general $0\le\al\le 2$.
The proof of Proposition \ref{prop2} is not difficult, essentially,
is based on a estimate of the following term
 \be
R_j\triangleq(S_{j-1}v\cdot\nabla)\tr_ju-\tr_j((v\cdot\nabla)u),
 \ee
 we give it in the following Lemma \ref{lem2.5}.
\end{rem}

\begin{lem}\label{lem2.5}
We rewrite $R_j$ as
$R_j=(S_{j-1}v-v)\cdot\nabla\tr_ju-[\tr_j,v\cdot\nabla]u$. Under the
condition (\ref{Condition}), there exists a sequence $c_j\in
l^q(\mb{Z})$ satisfying $\|c_j\|_{l^q}=1$, such that
 \bes
2^{js}\|R_j\|_{L^p}\le C\|\nabla v\|_{\dot
B^{n/p_1}_{p_1,\infty}\cap L^\infty}\|u\|_{\dot B^s_{p,q}}
 \ees
for any $j\in \mb{Z}$, where $C=C(n,q,s,p,p_1)$ is a constant
depending only on $n,q,s,p,p_1$.
\end{lem}

R. Danchin in \cite{Danchin} proved Lemma \ref{lem2.5} in the
inhomogeneous Besov space case. Similarly, using Bony's para-product
decomposition and Bernstein inequality, it is not difficult to prove
Lemma \ref{lem2.5}.

%%%%%%%%%%%%%%%%%%%%%%%%%%%%%%%%%%%%%%%%%%% section 3

\section{Local well-posedness and some blow-up criteria}
\setcounter{equation}{0}

In this section we prove Theorem 1.2 and Theorem 1.3, which are the
local well-posedness and the blow-up criteria of smooth solutions.

\vskip0.1in {\bf Step 1} Linear approximate equations

We construct sequence of approximate solutions by the following
linear equations:
 \be
\begin{cases}
\frac{\pt\theta^{k+1}}{\pt
t}+u^k\cdot\nabla\theta^{k+1}+\nu\La\theta^{k+1}=0,\ x\in\mr^3,\
t>0,\\
u^k=c(\theta^k)+\mc{P}(\theta^k), \ x\in\mr^3,\
t>0,\\
\mx{div }u^{k+1}=0,\\
\theta^{k+1}(0,x)=\theta_0(x),\ x\in\mr^3.
\end{cases}
 \ee
We set $\theta^0\triangleq \mx{e}^{-\nu t\La}\theta_0$, obviously,
$\theta^0\in L^1(\mr^+;\dot B^{3/p+1}_{p,1}(\mr^3))$. By
Proposition \ref{prop2}, we thus have
 \be
\theta^k\in \tilde L^\infty(\mr^+;\dot B^{3/p}_{p,1}(\mr^3))\cap
L^1(\mr^+;\dot B^{3/p+1}_{p,1}(\mr^3)),
 \ee
for any $k\ge 1$.

\vskip0.1in {\bf Step 2} Uniform estimates

We also need to obtain a uniform bound of $\theta^k(t,x)$ in
$\tilde L^\infty(\mr^+;\dot B^{3/p}_{p,1}(\mr^3))\cap
L^1(\mr^+;\dot B^{3/p+1}_{p,1}(\mr^3))$ for some a $T>0$
independent on $k$.

By the standard local existence method it is not difficult to prove
that there exists some time $T$ depended on the profile of
$\theta_0$ such that
 \be
\int^T_0\|\theta^k(t)\|_{\dot B^{3/p+1}_{p,1}}\md t\le C_0,
 \ee
for $k\ge 1$. For details refer to \cite{A-H,M-W}.

 In proposition \ref{prop2}, we take $r=1$ and $r=\infty$,
respectively, noting the Sobolev embedding relation $\dot
B^{3/p}_{p,1}(\mr^3)\hookrightarrow L^\infty(\mr^3)$ and the
boundedness of singular integral operator $\mc{P}$ on homogeneous
Besov space $\dot B^{3/p}_{p,1}(\mr^3)$ , it yields
 \be\nonumber
\|\theta^{k+1}\|_{\tilde L^\infty_T\dot
B^{3/p}_{p,1}}+\|\theta^{k+1}\|_{L^1_T\dot B^{3/p+1}_{p,1}}&\le&
C\exp\Bigg\{C\int^T_0\|u^k(\tau)\|_{\dot
B^{3/p+1}_{p,1}}\md\tau\Bigg\}\|\theta_0\|_{\dot B^{3/p}_{p,1}}\\
&\le& C \|\theta_0\|_{\dot B^{3/p}_{p,1}}.
 \ee

Consequently, the sequence $\{\theta^k\}$, $k\in \mb{N}$ is
uniformly bounded in $\tilde L^\infty(\mr^+;\dot
B^{3/p}_{p,1}(\mr^3))\cap L^1(\mr^+;\dot B^{3/p+1}_{p,1}(\mr^3))$.

\vskip0.1in {\bf Step 3} Strong convergence

We prove that $\{\theta^k\}$, $k\in \mb{N}$ is a Cauchy sequence
in $\tilde L^\infty(\mr^+;\dot B^{3/p}_{p,1}(\mr^3))\cap
L^1(\mr^+;\dot B^{3/p+1}_{p,1}(\mr^3))$, so it is strong
convergent.

Let $n,\ m\in \mb{N}$, and $n>m$. Set
$\theta^{n,m}\triangleq\theta^n-\theta^m$ and $u^{n,m}\triangleq
u^n-u^m= c(\theta^n-\theta^m)+\mc{P}(\theta^n-\theta^m)$. A simple
deduction yields
 \be
\begin{cases}
\pt_t\theta^{n+1,m+1}+u^n\cdot\nabla\theta^{n+1,m+1}+\nu\La\theta^{n+1,m+1}
=-u^{n,m}\cdot\nabla\theta^{m+1},\ t>0,\ x\in\mr^3\\
\mx{div}u^n=\mx{div}u^m=0,\\
\theta^{n+1,m+1}(0,x)=0,\ x\in\mr^3.
\end{cases}
 \ee

According to Proposition \ref{prop2}, noting the Sobolev embedding
relation $\dot B^{3/p}_{p,1}(\mr^3)\hookrightarrow
L^\infty(\mr^3)$ and the boundedness of singular integral operator
$\mc{P}$ on homogeneous Besov space $\dot B^{3/p}_{p,1}(\mr^3)$,
one has
 \be\nonumber
\|\theta^{n+1,m+1}\|_{\tilde L^\infty_T\dot B^{3/p}_{p,1}}\le
C\exp\Bigg\{C\int^T_0\|u^n(\tau)\|_{\dot
B^{3/p+1}_{p,1}}\md\tau\Bigg\}\int^T_0\|u^{n,m}\cdot\nabla\theta^{m+1}(
\tau )\|_{\dot B^{3/p}_{p,1}}\md\tau.\\
\label{3.6}
 \ee
By (\ref{product1}) in Proposition \ref{prop1} and the Sobolev
embedding relation $\dot B^{3/p}_{p,1}(\mr^3)\hookrightarrow
L^\infty(\mr^3)$, it can be deduced
 \bes
\|u^{n,m}\cdot\nabla\theta^{m+1}(\md \tau )\|_{\dot
B^{3/p}_{p,1}}\le C\|u^{n,m}\|_{\dot
B^{3/p}_{p,1}}\|\theta^{m+1}\|_{\dot B^{3/p+1}_{p,1}}.
 \ees
By the boundedness of singular integral operator $\mc{P}$ on
homogeneous Besov space $\dot B^{3/p}_{p,1}(\mr^3)$, we have
 \be\label{3.8}
\|u^{n,m}\|_{\dot B^{3/p}_{p,1}}\le C\|\theta^{n,m}\|_{\dot
B^{3/p}_{p,1}},\ \mx{ and } \|u^n(\tau)\|_{\dot
B^{3/p+1}_{p,1}}\le C\|\theta^n(\tau)\|_{\dot B^{3/p+1}_{p,1}}.
 \ee
Substituting (\ref{3.8}) into (\ref{3.6}), and choosing $T$ small
enough if necessary, we arrive at
 \bes
\|\theta^{n+1,m+1}\|_{\tilde L^\infty_T\dot B^{3/p}_{p,1}}&\le&
C\|\theta^{n,m}\|_{\tilde L^\infty_T\dot
B^{3/p}_{p,1}}\mx{e}^{\|\theta^n\|_{L^1_T\dot
B^{3/p+1}_{p,1}}}\int^T_0\|\theta^{m+1}\|_{\dot B^{3/p+1}_{p,1}}\md t\\
&\le& \ep \|\theta^{n,m}\|_{\tilde L^\infty_T\dot B^{3/p}_{p,1}},
 \ees
with $\ep<1$. Arguments by induction yield
 \be\label{3.81}
\|\theta^{n+1,m+1}\|_{\tilde L^\infty_T\dot B^{3/p}_{p,1}}\le
\ep^{m+1}\|\theta^{n-m,0}\|_{\tilde L^\infty_T\dot
B^{3/p}_{p,1}}\le C\ep^{m+1}\|\theta_0\|_{\dot B^{3/p}_{p,1}}.
 \ee
 Arguing similarly to above by Proposition \ref{prop2} it can be
 derived that
 \be\label{3.9}
\|\theta^{n+1,m+1}\|_{L^1_T\dot B^{3/p+1}_{p,1}}\le
\ep^{m+1}\|\theta^{n-m,0}\|_{\tilde L^\infty_T\dot B^{3/p}_{p,1}}\le
C\ep^{m+1}\|\theta_0\|_{\dot B^{3/p}_{p,1}}.
 \ee
Estimates (\ref{3.81})-(\ref{3.9}) imply that $\{\theta^k\}$ is a
Cauchy sequence in $\tilde L^\infty([0,T);\dot
B^{3/p}_{p,1}(\mr^3))\cap L^1((0,T);\dot B^{3/p+1}_{p,1}(\mr^3))$
for $k=0,1,\cdots$. Thus there exists a $\theta\in \tilde
L^\infty([0,T);\dot B^{3/p}_{p,1}(\mr^3))\cap L^1((0,T);\dot
B^{3/p+1}_{p,1}(\mr^3))$ such that $\theta^k$ converges strongly
to $\theta$.

\vskip0.1in {\bf Step 4} Uniqueness

 Suppose $\theta_1$ and
$\theta_2$ are two solutions of the equations (\ref{PE}) with the
same initial data $\theta_0\in \dot B^{3/p}_{p,1}$, and $\theta_1,\
\theta_2\in \tilde L^\infty([0,T);\dot B^{3/p}_{p,1}(\mr^3))\cap
L^1([0,T);\dot B^{3/p+1}_{p,1}(\mr^3))$. Introducing notations
$\theta_{1,2}\triangleq \theta_1-\theta_2$ and $u_{1,2}\triangleq
u_1-u_2=c(\theta_1-\theta_2)+\mc{P}(\theta_1)-\mc{P}(\theta_2)$ one
has by a simple deduction
 \be
\begin{cases}
\pt_t\theta_{1,2}+u_1\cdot\nabla\theta_{1,2}+\nu\La\theta_{1,2}
=-u_{1,2}\cdot\nabla\theta_2,\\
\mx{div} u_1=\mx{div}u_2=0,\\
\theta_{1,2}(0,x)=0.
\end{cases}
 \ee
Arguing similarly to the above (\ref{3.6}) we obtain
 \be
\|\theta_{1,2}\|_{\tilde L^\infty_t\dot B^{3/p}_{p,1}}\le
C\mx{e}^{C\|\theta_1\|_{L^1_t}\dot
B^{3/p+1}_{p,1}}\int^t_0\|\theta_{1,2}\|_{\tilde L^\infty_\tau\dot
B^{3/p}_{p,1}}\|\theta_2\|_{\dot B^{3/p+1}_{p,1}}\md\tau.
 \ee
Gronwall's inequality implies that $\theta_1(t)=\theta_2(t)$ for
any $0\le t\le T$.

\vskip0.1in {\bf Step 5} Smoothing effect

We shall prove the following regularity estimate
 \be\label{3.11}
\|t^\beta\theta\|_{\tilde L^\infty_T\dot B^{3/p+\beta}_{p,1}}\le
C(\beta)\mx{e}^{C\beta\|\theta\|_{L^1_T\dot
B^{3/p+1}_{p,1}}}\|\theta\|_{\tilde L^\infty_T\dot B^{3/p}_{p,1}},
 \ee
 for $\beta\ge 0$, where $t^\beta \theta$ obviously satisfies the following equation
 \be
\begin{cases}
\pt_t(t^\beta\theta)+(u\cdot\nabla)(t^\beta\theta)+\nu\La(t^\beta\theta)=-\beta
t^{\beta-1}\theta,\ t>0,\ x\in \mr^3,\\
\mx{div}u=0,\\
 (t^\beta\theta)(0,x)=0,\ x\in\mr^3.
\end{cases}
 \ee
We prove the estimate (\ref{3.11}) by induction. When $\beta=1$,
Proposition \ref{prop2} implies
 \be
\|t\theta\|_{\tilde L^\infty_T\dot B^{3/p+1}_{p,1}}\le
C\mx{e}^{C\|\theta\|_{L^1_T\dot B^{3/p+1}_{p,1}}}\|\theta\|_{\tilde
L^\infty_T\dot B^{3/p}_{p,1}}.
 \ee
Here Sobolev imbedding relation $\dot
B^{3/p}_{p,1}(\mr^3)\hookrightarrow L^\infty(\mr^3)$ and the
boundedness of singular integral operator $\mc{P}$ on homogeneous
Besov space $\dot B^{3/p+1}_{p,1}(\mr^3)$ have been used.

Assume the estimate (\ref{3.11}) is true for $\beta=k$, then
Proposition \ref{prop2} and induction implies, for $\beta=k+1$,
that
 \be\nonumber
\|t^{k+1}\theta\|_{\tilde L^\infty_T\dot B^{3/p+k+1}_{p,1}}&\le&
C\mx{e}^{C\|\theta\|_{L^1_T\dot
B^{3/p+1}_{p,1}}}\|t^k\theta\|_{\tilde L^\infty_T\dot
B^{3/p+k}_{p,1}}\\&\le& C(k+1)\mx{e}^{C(k+1)\|\theta\|_{L^1_T\dot
B^{3/p+1}_{p,1}}}\|\theta\|_{\tilde L^\infty_T\dot B^{3/p}_{p,1}}.
 \ee
For general $\beta\ge0$, noting $[\beta]\le \beta\le [\beta]+1$,
by interpolation between $\tilde L^\infty_T\dot
B^{3/p+[\beta]}_{p,1}(\mr^3)$ and  $\tilde L^\infty_T\dot
B^{3/p+[\beta]+1}_{p,1}(\mr^3)$,
 \be
\|t^\beta\theta\|_{\tilde L^\infty_T\dot B^{3/p+\beta}_{p,1}}\le
C\|t^{[\beta]}\theta\|^{[\beta]+1-\beta}_{\tilde L^\infty_T\dot
B^{3/p+[\beta]}_{p,1}}\|t^{[\beta]+1}\theta\|^{\beta-[\beta]}_{\tilde
L^\infty_T\dot B^{3/p+[\beta]+1}_{p,1}},
 \ee

we immediately prove the estimate (\ref{3.11}) for any
$\beta\ge0$.

\vskip0.1in {\bf Proof of Theorem \ref{thm1.3}}

In the proof of uniform estimate of the approximation solution
sequence $\theta^n$, we have obtained that if
 \be
\sum_{j\in\mb{Z}}(1-\mx{e}^{-C(T-\tau)2^j})^{1/2}\|\tr_j\theta(\tau)\|_{L^\infty}\le
\ep_0,
 \ee
for some a constant $\ep_0$. Then the solution $\theta$ is uniform
bounded
 \be
\|\theta\|_{\tilde L^2_T\dot
B^{3/p+1/2}_{p,1}}+\|\theta\|_{L^1_T\dot B^{3/p+1}_{p,1}}\le 2\ep_0,
 \ee
the solution $\theta$ thus can be extended beyond $T$. For details
see \cite{A-H,M-W}.

Let $[0,T^*)$ be the maximal existence interval, if $T^*<\infty$,
then
 \be
\lim \inf_{\tau\rightarrow
T^*}\sum_{j\in\mb{Z}}(1-\mx{e}^{-C(T-\tau)2^j})^{1/2}\|\tr_j\theta(\tau)\|_{L^\infty}
\ge\ep_0,
 \ee
otherwise, it can be extended beyond $T^*$. Noticing
$\|\tr_j\theta\|_{L^\infty}\le C\|\theta_0\|_{L^\infty}$, by the
Bernstein inequality (\ref{2.3}) we have
 \bes
\ep_0&\le& \lim \inf_{\tau\rightarrow T^*}\Big(\sum_{j\le
N}(1-\mx{e}^{-C(T^*-\tau)2^j})^{1/2}\|\tr_j\theta\|_{L^\infty}+\sum_{j\ge
N}(1-\mx{e}^{-C(T^*-\tau)2^j})^{1/2}\|\tr_j\theta\|_{L^\infty}\Big)\\
&\le& \lim \inf_{\tau\rightarrow
T^*}\Big((T^*-\tau)^{1/2}\|\theta_0\|_{L^\infty}\sum_{j\le
N}2^{j/2}+\|\nabla \theta\|_{L^\infty}\sum_{j\ge N}2^{-j}\Big)\\
&\le&\lim \inf_{\tau\rightarrow
T^*}\Big((T^*-\tau)^{1/2}\|\theta_0\|_{L^\infty}2^{N/2}+\|\nabla
\theta\|_{L^\infty}2^{-N}\Big).
 \ees
If we choose appropriate $N$, it follows
 \be
\lim \inf_{\tau\rightarrow
T^*}(T^*-\tau)\|\nabla\theta(\tau)\|_{L^\infty}\ge \ep_0.
 \ee

 Therefore, if
 \be
\int^T_0\|\nabla\theta(t)\|_{L^\infty}<\infty,
 \ee
there exists some a $T'>T$ such that $\theta$ can be continually
extended to $[0,T')$.

%%%%%%%%%%%%%%%%%%%%%%%%%%%%%section 4
\section{Global well-posedness}
In this section by virtue of the method of modulus of continuity
\cite{K-N-V} we prove the global well-posedness of Theorem
\ref{thm1.1}. In this case the difficulty is to construct a special
modulus of continuity which the solution $\theta$ has. First we
define a modulus of continuity.

\begin{defin}
Let $\om(x):\ [0,+\infty)\rightarrow [0,+\infty)$ be an increasing
continuous concave function satisfying $\om(0)=0$. We call  a
function $f$ from $\mr^n$ to $\mr^m$ has modulus of continuity
$\om$, if
 \bes
|f(x)-f(y)|\le \om(|x-y|), \mx{ for any } x,\ y\in\mr^n.
 \ees
\end{defin}

We recall that Kiselev, Nazarov and Volberg in \cite{K-N-V} proved a
lemma that said the Riesz transform didn't violate the modulus of
continuity too much as:

\begin{prop}\label{prop3}
If the function $\theta$ has a modulus of continuity $\om$, then
$u=(-R_2\theta, R_1\theta)$ has modulus of continuity
 \be
\Om_1(\xi)=A\bigg(\int^\xi_0\frac{\om(s)}s\md
s+\xi\int^\infty_\xi\frac{\om(s)}{s^2}\md s\bigg)
 \ee
 with a universal constant $A>0$, where $R_j$ is the $j-$th Riesz transform.
\end{prop}

We also need to prove that the singular integral operators
(\ref{1.3}) which are equivalent to double Riesz transforms or
their combinations do not spoil modulus of continuity too much,
although they do not preserve a modulus of continuity, see the
following Lemma \ref{lem3.1}.

\begin{lem}\label{lem3.1}
If the function $\theta$ has a modulus of continuity $\om$, then
$v=(R_1R_3\theta, R_2R_3\theta, -R^2_1\theta-R^2_2\theta)$ has
modulus of continuity
 \be\label{4.9}
\Om(\xi)=C\bigg(\int^\xi_0\frac{\om(s)}s\log\Big(\frac{\mx{e}\xi}
s\Big)\md
s+\xi\int^\infty_\xi\frac{\om(s)}{s^2}\log\Big(\frac{\mx{e}s}\xi\Big)\md
s\bigg),
 \ee
where $C>0$ is a constant.
\end{lem}

The proof is very simple, it only need a direct computation. Indeed,
Let $\Om_1(\xi)=A\bigg(\int^\xi_0\frac{\om(s)}s\md
s+\xi\int^\infty_\xi\frac{\om(s)}{s^2}\md s\bigg)$, then
 \bes
\Om(\xi)=C\bigg(\int^\xi_0\frac{\Om_1(\eta)}\eta\md
\eta+\xi\int^\infty_\xi\frac{\Om_1(\eta)}{\eta^2}\md \eta\bigg).
 \ees
By Fubini theorem, exchanging the order of the two integrals can
yield the result (\ref{4.9}) easily.

According to the blow-up criterion in Theorem \ref{thm1.3}, we need
to give a bound of $\|\nabla\theta\|_{L^\infty}$. For this purpose,
we choose the modulus of continuity $\om$ satisfying
$\om'(0)<\infty$ and $\lim_{\xi\rightarrow 0^+}\om''(\xi)=-\infty$.
Thus by the definition of modulus of continuity, it is not difficult
to prove that
 \be\label{bound}
\|\nabla\theta\|_{L^\infty}\le \om'(0).
 \ee

Let $T^*$ be the maximal existence time of the solutions $\theta\in
\tilde L^\infty([0,T^*);\dot B^{3/p}_{p,1}(\mr^3))\cap
L_{loc}^1((0,T^*);\dot B^{3/p+1}_{p,1}(\mr^3))$ to (\ref{PE}). By
Theorem \ref{thm1.2}, there exists a $T_0>0$ such that
 \bes
t\|\nabla\theta\|_{L^\infty}\le C\|\theta_0\|_{\dot B^{3/p}_{p,1}},
\mx{ for any } t\in [0,T_0].
 \ees
Let $\lambda>0$ and $T_1\in (0,T_0)$, we define
 \be
I=\{T\in [T_1,T^*): \forall t\in [T_1,T],
|\theta(t,x)-\theta(t,y)|<\om_\lambda(|x-y|),\ \mx {for any
}x\not=y\},
 \ee
where $\om_{\lambda}(\xi)\triangleq\om(\lambda\xi)$.

By appropriately choosing $\lambda$, for instance, set
$\lambda=\frac
{\om^{-1}(3\|\om_0\|_{L^\infty})}{2\|\om_0\|_{L^\infty}}
\|\nabla\theta(T_1)\|_{L^\infty}$, we can prove that $T_1\in I$, for
detail see \cite{M-W}. Thus $I$ is an interval of the form
$[T_1,T_*)$, where $T_*$ is the maximal of $T\in I$. We discuss the
relations between $T_*$ and $T^*$ in three cases, respectively.

{\bf Case 1:} If $T_*=T^*$, then in light of the inequality
(\ref{bound}) and the blow-up criterion (\ref{1.5}) in Theorem
\ref{thm1.3} we have $T^*=\infty$.

{\bf Case 2:} If $T_*\in I$, It is not difficult to prove that there
exists a positive $\eta$ such that $T_*+\eta\in I$, which is a
contradiction to the fact that $T_*$ is the maximal of $T\in I$, for
detail see \cite{M-W}.

{\bf Case 3:} If $T_*\not\in I$. The continuity of $\theta$ in time
implies that there exist $x\not=y$ such that
 \be
\theta(T_*,x)-\theta(T_*,y)=\om_\lambda(\xi),
 \ee
where $\xi=|x-y|$.

We shall prove that it is not possible. Let $f(t)\triangleq
\theta(t,y)-\theta(t,x)$ for the above fixed $x,y$. Clearly $f(t)\le
f(T_*)$ for any $t\in [0,T_*]$ by the definition of $I$. On the
other hand, we shall prove that $f'(T_*)<0$, which is a
contradiction.

The idea of proof is from \cite{K-N-V}, the difficulty is to
construct a modulus of continuity, for convenience of reading we
give a sketch of the proof.

By the regularity of solutions the equations can be defined in the
classical mean,
 \be
f'(T_*)=u(T_*,x)\cdot\nabla\theta(T_*,x)-u(T_*,y)\cdot\nabla\theta(T_*,y)
+\nu\La\theta(T_*,x)-\nu\La\theta(T_*,y).
 \ee

A direct computation by derivative immediately yields (see
\cite{K-N-V})
 \be
u(T_*,x)\cdot\nabla\theta(T_*,x)-u(T_*,y)\cdot\nabla\theta(T_*,y)\le
C(\om_\lambda(\xi)+\Om_\lambda(\xi))\om'_\lambda(\xi).
 \ee

Noting that the dissipative term $\La\theta(x,t)$ can be written
as $\frac {\md}{\md s}P_s*\theta|_{s=0}$, where
 \bes
P_s(x)=\frac {s}{\pi^2(|x|^2+s^2)^{3/2}}
 \ees
is the three dimensional Poisson kernel in $\mr^3$. By a detail
deduction (use of the symmetry and monotonicity of the Poisson
kernel and some integral techniques, see \cite{K-N-V}) we have
 \be\nonumber
\nu\La\theta(T_*,x)-\nu\La\theta(T_*,y)&\le&
\frac\nu\pi\int^{\frac\xi2}_0\frac{\om_\lambda(\xi+2s)
+\om_{\lambda}(\xi-2s)-2\om_\lambda(\xi)}{s^2}\md s\\\nonumber
&&+\frac\nu\pi\int_{\frac\xi2}^\infty\frac{\om_\lambda(\xi+2s)
-\om_{\lambda}(\xi-2s)-2\om_\lambda(\xi)}{s^2}\md s\\&\triangleq&
\lambda J(\lambda \xi),
 \ee
where
 \be\label{4.25}
J(\xi)&=&\frac\nu\pi\int_0^\frac{\xi}{2}\frac{\om(\xi+2s)
-\om(\xi-2s)-2\om(\xi)}{s^2}\md
s\\&&+\frac\nu\pi\int_{\frac\xi2}^\infty\frac{\om(\xi+2s)
-\om(\xi-2s)-2\om(\xi)}{s^2}\md s.
 \ee
Thus we only need to prove
 \be\label{4.26}
f'(T_*)\le C\lambda[(\om+\Om)\om'+J](\lambda\xi)<0.
 \ee

For this purpose we choose the modulus of continuity $\om$ as
follows
 \be
\om(x)=\begin{cases} x-x^{3/2},\ \ &\mx{ if }0\le x\le \de,\\
\de-\de^{3/2}+\frac\ga 3\arctan\frac{1+\log \frac x\de}3-\frac\ga
3\arctan\frac13,\ \ &\mx{ if } \de\le x.
\end{cases}
 \ee
Its derivative is
 \be\label{4.28}
\om'(x)=\begin{cases} 1-\sqrt x,\ \ &\mx{ if }0\le x< \de,\\
\frac \ga{x[9+(1+\log\frac x\de)^2]},\ \ &\mx{ if } \de< x.
\end{cases}
 \ee
Here $\de>\ga>0$ are two small enough constants that will be
determined later. Obviously $\om$ is concave and satisfies
 \bes
\om'(0)<\infty \mx{ and } \lim_{x\rightarrow 0^+}\om''(x)=-\infty.
 \ees
In the following we prove the inequality (\ref{4.26}) in two cases.

{\bf Case 1: $0\le\xi\le\de$.}

Since $\om(s)\le s$ for all $0\le s \le\de$, we have
 \be\label{4.29}
\int^\xi_0\frac{\om(s)}s\log\frac{\mx{e}\xi}s\md s\le
\int^\xi_0\log\frac{\mx{e}\xi}s\md s\le 2\xi,
 \ee\label{4.30}
 \be\nonumber
\xi\int^\de_\xi\frac{\om(s)}{s^2}\log\frac{\mx{e}s}\xi\md s&\le&
\xi\int^\de_\xi\frac{1}{s}\log\frac{\mx{e}s}\xi\md s\\
&=&\frac12\xi\log\frac\de\xi\bigg(2+\log\frac\de\xi\bigg),
 \ee
and
 \be\nonumber
&&\xi\int_\de^\infty\frac{\om(s)}{s^2}\log\frac{\mx{e}s}\xi\md s\le
\xi\frac{\om(\de)}\de\log\frac{\mx{e}\de}\xi+\xi\int^\infty_\de\bigg(\frac{\om'(s)}s
\log\frac{\mx{e}s}\xi+\frac{\om(s)}{s^2}\bigg)\md s\\\nonumber &\le&
\xi\log\frac{\mx{e}\de}\xi+\frac{\xi\ga}{4\de}
\log\frac{\mx{e}\de}\xi+\xi\int^\infty_\de\frac{\om'(s)}s\md s\\
\label{4.31} &\le&
\xi\bigg(1+\log\frac{\mx{e}\de}\xi\bigg)\bigg(1+\frac\ga{2\de}\bigg).
 \ee
Collecting (\ref{4.29})-(\ref{4.31}), we get
 \be
\om(\xi)+\Om(\xi)\le 3\xi+\xi\bigg(2+\log\frac\de\xi\bigg)^2.
 \ee

Next, we estimate the negative part $J$, only use of the first
integral in (\ref{4.25}) is enough. The concavity of $\om$, Taylor
formula and monotonicity of $\om''$ on $[0,\xi]$ imply that
 \be
\frac\nu\pi\int^{\frac\xi2}_0\frac{\om(\xi+2s)+\om(\xi-2s)-2\om(\xi)}{s^2}\md
s\le \frac1\pi\xi\om''(\xi)=-\frac{3\nu} {4\pi}\xi\xi^{-1/2}.
 \ee
Obviously, if $\xi\in(0,\de]$ and $\de>0$ is small enough, one has
 \be
(\om(\xi)+\Om(\xi))\om'+J(\xi)\le\xi\bigg[3+\Big(2+\log\frac\de\xi\Big)^2
-\frac{3\nu}{4\pi}\frac1{\sqrt\xi}\bigg]<0.
 \ee

{\bf Case 2: $\xi\ge\de$.}

In this case we have $\om(s)\le s$ for $0\le s\le\de$ and
$\om(s)\le\om(\xi)$ for $\de\le s\le\xi$. Therefore,
 \be\nonumber
&&\int^\xi_0\frac{\om(s)}s\log\frac{\mx{e}\xi}s\md
s=\bigg(\int_0^\de+\int^\xi_\de\bigg)\frac{\om(s)}s\log\frac{\mx{e}\xi}s\md
s\\ \nonumber &\le& \int^\de_0\log\frac{\mx{e}\xi}s\md
s+\int^\xi_\de\frac{\om(\xi)}s\log\frac{\mx{e}\xi}s\md
s\\&\le&\nonumber
\de\bigg(2+\log\frac\xi\de\bigg)+\om(\xi)\bigg(\log\frac\xi\de\bigg)
\bigg(1+\log\frac\xi\de\bigg)\\ \label{4.37} &\le&
\om(\xi)\bigg[1+\Big(1+\log\frac\xi\de\Big)^2\bigg],
 \ee
where we use the fact $\de\le\om(\de)\le\om(\xi)$.

Arguing similarly to above it can be derived that
 \be\nonumber
\xi\int^\infty_\xi\frac{\om(s)}{s^2}\log\frac{\mx{e}s}\xi\md
s&=&\om(\xi)+\xi\int^\infty_\xi\frac1s\bigg(\om(s)\log\frac{\mx{e}s}\xi\bigg)'\md
s\\ \nonumber
&\le&2\om(\xi)+\xi\ga\int^\infty_\xi\frac1{s^2}\log\frac{\mx{e}s}\xi\md
s+\xi\int^\infty_\xi\frac{\om'(s)}s\md s\\\label{4.38} &\le&
2\om(\xi)+3\ga\le 5\om(\xi).
 \ee
Combining the estimates (\ref{4.37})-(\ref{4.38}) with $\om'(\xi)$
of (\ref{4.28}) we get
 \be\nonumber
(\om(\xi)+\Om(\xi))\om'(\xi)&\le&
C\om(\xi)\bigg[7+\Big(1+\log\frac\xi\de\Big)^2\bigg]\frac
\ga{\xi[9+(1+\log\frac\xi\de)^2]}\\ &\le&C\ga\frac{\om(\xi)}\xi.
 \ee

To complete the proof, we only need to estimate the second integral
in $J$. In case $\de\le\xi$, we have
 \be
\om(2\xi)\le \om(\xi)+\om'(\eta)\xi\le \om(\xi)+\frac\ga9\le
\frac{10}9\om(\xi).
 \ee
The concavity of $\om(x)$ implies $\om(2s+\xi)-\om(2s-\xi)\le
\om(2\xi)$ for all $\frac\xi2\le s$, thus it reaches
 \be
\frac\nu\pi\int^\infty_{\frac\xi2}\frac{\om(2s+\xi)-\om(2s-\xi)-2\om(\xi)}{s^2}\md
s\le -\frac{16\nu}{9\pi}\frac{\om(\xi)}\xi.
 \ee
It follows that
 \be
(\om(\xi)+\Om(\xi))\om'+J(\xi)\le
\frac{\om(\xi)}\xi\Big(C\ga-\frac{16\nu}{9\pi}\Big)<0,
 \ee
if we take $\ga<\min\{\frac{16}{9\pi C},\de\}$. The proof of Theorem
\ref{thm1.1} is thus completed.

%%%%%%%%%%%%%%%%%%%%%%%%%%%%%%%%%%%%  Acknowledgements
\textbf{Acknowledgements} This work was completed during his visit
Beijing Institute of Applied Physics and Computational Mathematics,
and authors would like to thank Prof. Changxing Miao for some
valuable discussion on this topic. The research of B Yuan was
partially supported by the National Natural Science Foundation of
China (No. 10771052).
% Natural Science Foundation of Henan Province
%(No. 0611055500) and Doctor Fund of Henan Polytechnic University(No.
%B2008-62).
%%%%%%%%%%%%%%%%%%%%%%%%%%%%%%%%%%%% References

\end{document}